\documentclass[english]{article}
\usepackage[T1]{fontenc}
\usepackage[latin1]{inputenc}
\usepackage{babel}

\makeatletter

\providecommand{\LyX}{L\kern-.1667em\lower.25em\hbox{Y}\kern-.125emX\@}

\usepackage{amsmath,amsfonts,amssymb,amsthm,latexsym,amscd}
 \long\def\unmarkedfootnote#1{{\long\def\@makefntext##1{##1}\footnotetext{#1}}}

\newtheorem{theorem}{Theorem}

\newtheorem{corollary}{Corollary}
\newtheorem{lemma}{Lemma}

\newtheorem{definition}{Definition}

\makeatother

\title{
Random Liouville functions and normal sets  }
\author{ Alexander Fish }

\begin{document}
\maketitle
\unmarkedfootnote{2000 Mathematics Subject Classification: 11N64,
05D10.}
 We define a random Liouville function \( \lambda_Q
\) which depends on a random set \( Q \) of primes and prove that
\( A_Q = \{ n \in \mathbb{N} | \lambda_Q(n) = -1 \} \) is normal
almost everywhere. This fact enables us to generate a family of
normal sets such that the equation \( xy =z \) is not solvable
inside them. Additionally we prove that  equations \( xy=z^2, x^2
+ y^2 = square, x^2 - y^2 = square \) are solvable in any normal
set and for any equation \( xy=cn^2 \) ( \( c > 1 \), is not a
square ) there exists a normal set \( A_c \) such that the
equation is not solvable inside \( A_c \).

\section{ Introduction }

     With the familiar notion of normal numbers in mind, we shall call
an infinite binary sequence  \textit{normal} if any binary word
\(\omega\) of length \(|\omega|\) occurs in the sequence with the
right frequency: \(2^{-|\omega|}\). We have the natural bijection
between infinite \(\{0,1\}\)-sequences \(\lambda\) and the subsets
of the natural numbers \(A_{\lambda} = \{i|\lambda_i = 1\}\).   We
now have
\begin{definition}
 A set \(B \subset  \mathbb{N}\)  is called  \textit{normal} if the corresponding
\(\{0,1\}\) sequence is normal.
\end{definition}
  In this note we shall be interested in normal sets and the possibility
of solving diophantine equations with integers from a given, but
arbitrary, normal set.  We expect that there are many diophantine
equations (or systems of equations) which, if they are solvable at
all in integers, are solvable with integers chosen from a given
normal set. We call such equations \textit{N-regular}, and we
denote by DSN the family of N-regular equations (or systems of
equations).

  An equation, or system of equations, is called "partition-regular", if
for any finite partition of the natural numbers, the system is
solvable within one of the cells of the partition.  One of the
earliest examples of a partition-regular equation is  Schur's
equation: \(x+y=z\).  It is not hard to see that Schur's equation
is also N-regular. Rado in \cite{rado} classified all systems of
linear diophantine equations that are partition regular. Rado's
theorem implies the familiar van der Waerden theorem on existence
of arbitrarily long monochromatic arithmetic progressions in any
finite coloring of the natural numbers.

    Using Furstenberg's theorem regarding Rado's systems in \cite{furst1}, one can
obtain the analogous result for N-regularity: namely, any Rado
system of linear equations is in DSN.

    From the foregoing, we have many linear equations in DSN. But little
is known in the non-linear case. For example, it is an open
question as to whether the Pythagorean equation \(x^2 + y^2 =
z^2\) is in DSN. The purpose of this note is to show that the
equation \(xy=z\) is not in DSN. This equation is called the
multiplicative Schur equation. It is an easy consequence of
Schur's additive theorem that his multiplicative equation is also
partition-regular. In fact in any finite partition of
\(\mathbb{N}\) one can find solutions to both the additive and the
multiplicative equations in the same cell (\cite{berg1}). Thus
partition regularity does not imply N-regularity. To show that
\(xy=z\) is not in DSN we will use a construction of random normal
sets, based on a variant of the Liouville function \(\lambda(n)\)
from number theory. Recall
\begin{definition}
Liouville's function \( \lambda : \mathbb{N} \rightarrow \{ -1, 1
\} \) is defined as follows:
\[
\lambda(p_1^{e_1}p_2^{e_2} \ldots p_k^{e_k}) = (-1)^{e_1+e_2+ \ldots +e_k}
\]
where \( p_1, \ldots , p_k \) are primes.
\end{definition}
It is a well known and very deep question whether the set \( A =
\{ n \in \mathbb{N} | \lambda(n) = -1 \} \) forms a normal set,
see \cite{liouv_1} and \cite{liouv_2}.  It seems that at present
we are far away from resolving this outstanding problem. But just
for clarity, if the answer for the question is positive, then the
aforementioned set \( A \) gives us an example of a normal set
with no solution to the equation \( xy=z \).
\newline
For the following we will use random Liouville's function \( \lambda_Q \) which
is defined by random choice of subset \( Q \) inside \(P\) (prime numbers) as follows
\[
\lambda_Q(p_1^{e_1}p_2^{e_2}\ldots p_k^{e_k}) =
\lambda_Q(p_1)^{e_1} \lambda_Q(p_2)^{e_2} \ldots
\lambda_Q(p_k)^{e_k}
\]
and
\begin{displaymath}
\lambda_Q(p) = \left\{
                       \begin{array}{cc} -1 &  p \in Q \\
                                   1 &  p \not \in Q
                                   \end{array}
                                   \right.
\end{displaymath}
 By randomness of \( Q \) we mean  that a choice of
every prime number \( p \)  is independent of other prime numbers
and \( Pr ( p \in Q ) = 0.5 \)  for any \( p \in P \).
\newline
One defines \(A_Q = \{ n \in \mathbb{N} | \lambda_Q(n) = -1 \}\).
In the second section of the note we prove  the following
\begin{theorem}
\label{main_theorem} For almost every \( Q \) the set \( A_Q \) is
normal.
\end{theorem}
This theorem gives us an infinite family of normal sets such that
the multiplicative Schur's equation is not solvable in these sets.
\newline
In the third part of the note we prove that the equations \(
xy=z^2 \), \( x^2 + y^2 = square\) and \( u^2 - v^2 = square \)
are in \(DSN \).

\textit{ Acknowledgment: } The author would like to thank Prof.
Hillel Furstenberg without whose support this work would not have
been done, Prof. Vitaly Bergelson for valuable discussions, Prof.
 Alex Samorodnitsky for remarkable suggestions, Michael Fish for help in proving lemma \ref{th2_pos} .
\section{ \( A_Q \) is normal for a.e. \( Q \)}
We start from an obvious claim about normality of \( A_Q \).
\begin{lemma}
\label{norm_lemma} Let \( Q \subset P \) be given, then \( A_Q \)
is a normal set \( \Leftrightarrow \) for any \(k \in (\mathbb{N}
\cup \{0\}) \) and any \( i_1 < i_2 <  \ldots < i_k \) we have
\begin{equation}
 \lim_{N \rightarrow \infty} \frac{1}{N} \sum_{n=1}^N \lambda_Q(n) \lambda_Q(n+i_1) \ldots \lambda_Q(n+i_k)
 =0.
 \label{eq:main_eq}
\end{equation}
\end{lemma}
We proceed with the following statement which is readily proved
\begin{lemma}
\label{tec_lemma} Let \( \{ a_n \} \) be a bounded sequence.
Denote by \( T_N = \frac{1}{N} \sum_{n=1}^N a_n \). Then \( T_N \)
converges to a limit \( t \) \( \Leftrightarrow \) there exists a
sequence of increasing indices \( \{ N_i \} \) such that \(
\frac{N_i}{N_{i+1}} \rightarrow 1 \) and \( T_{N_i} \rightarrow_{i
\rightarrow \infty} t \).
\end{lemma}
The next step is to show
\begin{displaymath}
 \sum_{N=1}^{\infty} E((\frac{1}{N^{40}} \sum_{n=1}^{N^{40}} \lambda_Q(n) \lambda_Q(n+i_1) \ldots \lambda_Q(n+i_k))^2) <
 \infty.
\end{displaymath}
\begin{lemma}
\label{main_lemma}
 Let \( T_N  \) as was defined previously, then
\( E( T_N^2 ) \leq O(\frac{1}{N^{0.05}}) \).
\end{lemma}
\begin{proof}
By linearity of expectation we get
\[
 E(T_N^2) = \frac{1}{N^2} \sum_{x,y=1}^{N} E(\lambda_Q(x) \lambda_Q(x+i_1) \ldots \lambda_Q(x+i_k)
\lambda_Q(y) \lambda_Q(y+i_1) \ldots \lambda_Q(y+i_k)).
\]
Note that for any \( m \in \mathbb{N} \), \( E(\lambda_Q(m)) = 0 \) unless \( m \) is a square in which case
\( E(\lambda_Q(m)) = 1 \).
\newline
Let us denote by
\[ \phi(x) \risingdotseq \lambda_Q(x) \lambda_Q(x+i_1) \ldots
\lambda_Q(x+i_k) \]
 and
 \[ \xi(x) \risingdotseq x(x+i_1) \ldots
(x+i_k). \]
 By distribution of \( Q \) we get
 \[ E (\phi(x)
\phi(y)) = 1    \Leftrightarrow  \xi(x) \xi(y) = m^2. \]
 Otherwise
\[ E (\phi(x) \phi(y)) = 0. \]
 Therefore, to obtain  an upper bound
on \( E({T_N}^2) \), we give an upper bound on the number of pairs
\( (x,y) \in [1,N] \times [1,N] \) which satisfy \( \xi (x) \xi(y)
= square \).
\newline
For a given \( x \in [1,N] \) let us assume that \( \xi(x) = c_x
m^2 \), where \( c_x \) is a square-free number, say \( c_x =
p_{j_1} \ldots p_{j_l} \) is the prime factorization of \( c_x \).
Then we will define \( h(x) = l \) (thus \( h(x) \) is a number of
primes in prime factorization of maximal square-free number which
divides \( x \)). Denote by \( D \) the set of all possible common
divisors  of the numbers \( x, x+i_1, \ldots, x+i_k \) (i.e.
positive integers which divide at least two of them).  For a
finite non empty set \( S \) of positive numbers we denote by \(
m(S) \) the product of all elements of \( S \)  and, for empty
set, we fix  \( m(\emptyset) = 1 \).
\newline
 Note that \( \xi(x) \xi(y) = square \) \( \Rightarrow \)
there exist \( S_1 \subset D \) and \( S_2 \subset \{p_{j_1},
\ldots ,p_{j_l} \}\) such that \( y = m(S_1) m(S_2) square \).
\newline
Assume  \( |D| = r \) (\( r \) depends only on the set \( \{i_1,
\ldots, i_k \} \) and doesn't depend on \(x\)). Then we obtain \(
\xi(x) \xi(y) = square \) for at most \( 2^r2^{h(x)} \sqrt{N} \)
\( y \)'s inside \( [1,N] \).
 Thus
\[
E(T_N^2) \leq \frac{1}{N^2} ( \sum_{n=1}^N 2^r2^{h(n)} \sqrt{N} )
\leq  \frac{c}{N^{1.5}} \sum_{n=1}^N 2^{h(n)}
\]
Therefore it remains to bound the expression \( \sum_{n=1}^N 2^{h(n)} \).
\newline
Let \( p = p_i \) be a smallest prime number such that \(
\frac{k+1}{\log_2{p}} \leq 0.45 \).
 If \( \xi(n) \) does not contain as dividers \( 2,3, \ldots, p \) then
\( h(n) \leq \log_p{(n+i_k)^{k+1}} = (k+1)
\frac{\log_2{(n+i_k)}}{\log_2{p}} \). This gives us
\[
2^{h(n)} \leq  (n+i_k)^{\frac{k+1}{\log_2{p}}} \leq (n+i_k)^{0.45}
\]
But if \( \xi(n) \) is arbitrary then \( h(n) \) can increase by
at most \( i \), this means \( 2^{h(n)} \leq 2^i (n+i_k)^{0.45}
\). Thus \( \sum_{n=1}^N 2^{h(n)} \leq C_1 (N+i_k)^{1.45} \) and
therefore we get
\[
E(T_N^2) \leq C_2 \frac{1}{N^{0.05}}
\]
\end{proof}
\begin{proof}(\textbf{theorem \ref{main_theorem})}
From the last lemma we conclude that \(  \sum_{N=1}^{\infty}
E(T_{N^{40}}^2) < \infty \). Thus almost surely \( T_{N^{40}}
\rightarrow 0\). By lemma \( \ref{tec_lemma} \) follows that
almost surely \( T_N \rightarrow 0\). And from lemma \(
\ref{norm_lemma} \) (and countability of necessary conditions) it
follows that for almost all \( Q \subset P \) the sets \( A_Q \)
are normal.
\end{proof}
We can now demonstrate the main result of this note.
\begin{theorem}
There exists \( A \subset \mathbb{N} \)  a normal set such that
the multiplicative Schur's equation is not solvable inside \( A
\).
\end{theorem}
\begin{proof}
 We have already shown the existence of many \(Q \) ( \( Q \subset P \))  such that \( A_Q \) are normal.
By definition of \( A_Q \) follows that for any \( x,y \in A_Q \)
the number \( xy \not \in A_Q \). Therefore we can't find \( x,y,z
\in A_Q \) such that \( xy = z \).
\end{proof}
\begin{corollary}
For any equation \( xy = c n^k \) (where \( c,k \) are natural
numbers, \( c \) is not a square and \( k \) is even) we can find
a normal set \( A_{c,k} \subset \mathbb{N} \) such that for any \(
x,y \in A \) we have \( xy \not = c n^k \) for every natural \( n
\).
\end{corollary}
\begin{proof}
We take \( A_Q \)  be a normal and such that \( \lambda_Q(c) = -1
\) (it happens with the positive probability \(\frac{1}{2} \), and
thus there exist such sets). Then obviously we can't solve the
proposed equation inside \( A_Q \).
\end{proof}

\section{ Solvability of equation \( xy = z^2 \) and related problems}
\begin{theorem}
\label{th1_pos} Let \( A \subset \mathbb{N} \) be a normal set.
Then there exist \( x,y,z \in A \) (\( x \not = y \)) such that \(
xy = z^2 \).
\end{theorem}
\begin{proof}
For a set \( S \subset \mathbb{N} \) let us define \( S_a = \{ n
\in \mathbb{N} | an \in S \} \), where \( a \in \mathbb{N} \). It
is easily seen that if \( S \) is normal then \( S_a \) is normal
for any natural \( a \) (see \cite{furst2}). We denote by \( d(S)
\) the density of a set \( S \), if it exists.
\newline
Let \( A \) be a normal set. We denote by \( R_n \risingdotseq
A_{2^n} \). For any \( n \) holds \( d(R_n) = \frac{1}{2} \). Let
us denote by
\[\mu_N (S) = \frac{|S \cap \{ 1, 2, , \ldots N
\}|}{N} \]
 for any \( S \subset \mathbb{N} \) and any \( N \in
\mathbb{N} \).
\newline
By Szemer\'{e}di's theorem (finite version), for any \( \delta >
0\) and any \( l \in \mathbb{N} \) there exists \( N(l,\delta) \)
such that for any \( N \geq N(l,\delta) \) and any \( F \subset
\{1,2,\ldots,N\} \) such that \( \frac{|F|}{N} \geq \delta \) the
set \( F \) contains an arithmetic progression of length \( l \)
(see \cite{szemeredi}).
\newline
One chooses \( K \geq  N(3,\frac{1}{3}) \). Then there exists \(
N_K \) such that \( \mu_{N_K} (R_i) \geq \frac{1}{3} \) for every
\( 1 \leq i \leq K \).
\newline
We claim that there exists \( F \subset \{1,2, \ldots, K \} \)
such that \( \frac{|F|}{K} \geq \frac{1}{3} \) and \(
\mu_{N_K}(\cap_{j \in F} R_j) > 0 \). If not, let us denote \(
1_{R_i} \) to be the indicator function of the set \( R_i \)
inside the set \( \{ 1, \ldots , N_K \} \). Then \[ \int_{[1,N_K]}
(1_{R_1} + \ldots + 1_{R_K}) d\mu_{N_K}  = \sum_{j=1}^{K}
\int_{[1,N_K]} 1_{R_j} d\mu_{N_K}  \geq \frac{K}{3} \]
 on the one hand.
\newline
But on the other hand
\[ \int_{[1,N_K]} (1_{R_1} + \ldots +
1_{R_K}) d\mu_{N_K}  < \frac{K}{3} \]
 because the function \(
1_{R_1} + \ldots + 1_{R_K} < \frac{K}{3} \).
\newline
Let \( F \subset \{1,2, \ldots, K \} \) such that \( \frac{|F|}{K} \geq \frac{1}{3} \) and
\( \mu_{N_K}(\cap_{j \in F} R_j) > 0 \). Then by the choice of \( K \) follows that \( F \) necessary contains arithmetic
progression of length \( 3 \). The last means there exist \( a,b,c \in F \) such that \( a+c=2b \). Let us take
\( R_a, R_b , R_c \). We have \( R_a \cap R_b \cap R_c \neq \emptyset \) and this means there exists \( n \in \mathbb{N} \)
such that \( n2^a \in A \) and \( n2^b \in A \) and \( n2^c \in A \). Let us denote by \( x,y,z \) the following elements of \( A \):
 \( x = n2^a \), \( y = n2^c \), \( z = n2^b \). Then we have
\[
xy = z^2
\]
\end{proof}
\textbf{Question:} Are the equations \( x y = c^2 z^2 \), where \(
c > 0 \) is a natural number, always solvable inside an arbitrary
normal set?
\newline
\begin{theorem}
\label{th2_pos}
 Let \( A \subset \mathbb{N} \) be an arbitrary
normal set. Then there exist \( x,y,u,v \in A \) such that \( x^2
+ y^2 = square \) and \( u^2 - v^2 = square \).
\end{theorem}
\begin{proof}
Note that there exist \( a,b,c \in \mathbb{N} \) such that \( a^2 + b^2 = square \) and
\( a^2 + c^2 = square \) and \( b^2 + c^2 = square \). For example \( a = 44, b = 117, c = 240 \).
\newline
Let \( A \subset \mathbb{N} \) be an arbitrary normal set. We look
at \( A_a, A_b, A_c \) which are defined as in \ref{th1_pos}. Then
\( d(A_a) = d(A_b) = d(A_c) = \frac{1}{2} \) and thus it can not
be true that the intersection of each pair from the triple is
empty.
\newline
Without loss of generality, let us assume that \( A_a \cap A_b \neq \emptyset \).
\newline
Thus there exists \( z \in A_a \cap A_b \) or equivalently \( z a, z b  \in A \). But \( a^2 + b^2 = square \) and therefore
\( (za)^2 + (zb)^2 = square \).
\newline
A proof that the equation \(  u^2 - v^2 = square \) is solvable in
any normal set is similar. We use the fact that there exist \( a,
b ,c \in \mathbb{N} \) such that \( a < b < c \) and holds \( c^2
- b^2 = square \) and \( c^2 - a^2 = square \) and
 \(b^2 - a^2 = square \). For example \( a = 153, b = 185, c = 697 \).
\end{proof}
\textbf{Questions:} 1) For an arbitrary normal set \( A \) do
there exist \( x,y,z \in A \) such that \( x^2 + y^2 = z^2\)?
\newline
2)  For an arbitrary normal set \( A \) do there exist \( x,y,z
\in A \) such that \( x^2 - y^2 = z^2\)?

\textit{Einstein Institute of Mathematics \\ The Hebrew University
of Jerusalem \\ Jerusalem, 91904, Israel \\
E-mail: afish@math.huji.ac.il}


\begin{thebibliography}{150}
\bibitem{berg1} Bergelson, V.; Hindman, N. Additive and multiplicative
Ramsey theorems in $N$---some elementary results. Combin. Probab.
Comput. 2 (1993), no. 3, 221--241.
 \bibitem{liouv_1} Cassaigne, J.; Ferenczi, S.; Mauduit, C.; Rivat, J.; Sßrk\"zy, A.
 On finite pseudorandom binary sequences. III. The Liouville function. I. Acta Arith. 87 (1999), no. 4, 367--390.
 \bibitem{liouv_2} Cassaigne, J.; Ferenczi, S.; Mauduit, C.; Rivat, J.; Sßrk\"zy, A.
  On finite pseudorandom binary sequences. IV. The Liouville function. II. Acta Arith. 95 (2000), no. 4, 343--359.
\bibitem{furst1} Furstenberg, H. Recurrence in Ergodic Theory and Combinatorial Number Theory. Princeton University Press, 1981.
\bibitem{furst2}  Furstenberg, H. Disjointness in ergodic theory, minimal sets, and a problem in Diophantine approximation.
 Math. Systems Theory 1 (1967), 1--49.
\bibitem{rado} Rado, R. Note on combinatorial analysis. Proc. London Math. Soc. 48 (1943), 122--160.
\bibitem{szemeredi} Szemer\'{e}di, E. On sets of integers containing no $k$ elements in arithmetic progression.
 Collection of articles in memory of Juri\v{i} Vladimirovi\v{c} Linnik. Acta Arith. 27 (1975), 199--245.
 \end{thebibliography}
\end{document}